\documentclass[12pt,leqno,draft]{article}

\def\dist{\mathop{\rm dist}}

\newtheorem{theorem}{Theorem}
\newtheorem{lemma}[theorem]{Lemma}
\newtheorem{proposition}[theorem]{Proposition}
\newtheorem{sublemma}[theorem]{Sublemma}
\newtheorem{definition}[theorem]{Definition}
\newtheorem{corollary}[theorem]{Corollary}
\newtheorem{problem}[theorem]{Problem}
\newtheorem{remark}[theorem]{Remark}
\newtheorem{claim}[theorem]{Claim}
\newtheorem{assumptions}[theorem]{Assumptions}
\newtheorem{examples}[theorem]{Examples}
\newtheorem{question}[theorem]{Question}
\newtheorem{sassumptions}[theorem]{Standing Assumptions}
\newtheorem{sassumption}[theorem]{Standing Assumption}
\newtheorem{conjecture}[theorem]{Conjecture}

\newcommand{\begintheorem}{\addtocounter{equation}{1}\begin{theorem}}
\newcommand{\beginlemma}{\addtocounter{equation}{1}\begin{lemma}}
\newcommand{\beginproposition}{\addtocounter{equation}{1}\begin{proposition}}
\newcommand{\beginsublemma}{\addtocounter{equation}{1}\begin{sublemma}}
\newcommand{\begindefinition}{\addtocounter{equation}{1}\begin{definition}}
\newcommand{\begincorollary}{\addtocounter{equation}{1}\begin{corollary}}
\newcommand{\beginproblem}{\addtocounter{equation}{1}\begin{problem}}
\newcommand{\beginremark}{\addtocounter{equation}{1}\begin{remark}}
\newcommand{\beginclaim}{\addtocounter{equation}{1}\begin{claim}}
\newcommand{\beginassumptions}{\addtocounter{equation}{1}\begin{assumptions}}
\newcommand{\beginexamples}{\addtocounter{equation}{1}\begin{examples}}
\newcommand{\beginquestion}{\addtocounter{equation}{1}\begin{question}}
\newcommand{\beginsassumptions}{\addtocounter{equation}{1}\begin{sassumptions}}
\newcommand{\beginsassumption}{\addtocounter{equation}{1}\begin{sassumption}}
\newcommand{\beginconjecture}{\addtocounter{equation}{1}\begin{conjecture}}

\begin{document}

\title{Some Remarks Concerning Integrals of Curvature on Curves and
Surfaces}

\author{S. Semmes\thanks{The author was partially supported by the
U.S. National Science Foundation.  A presentation based on this paper
was made at the AMS Special Session ``Surface Geometry and Shape
Perception'' (Hoboken, 2001).}}

\date{}

\maketitle

\begin{abstract}
In this paper we discuss some topics that came up in Chapters 2 and 3
of Part III of \cite{DS2}.  These involve relations between derivatives
of Cauchy integrals on curves and surfaces and curvatures of the
curves and surfaces.  In ${\bf R}^n$ for $n > 2$, ``Cauchy integrals''
can be based on generalizations of complex analysis using quarternions
or Clifford algebras (as in \cite{BDS}).  Part of the point here is to
bring out the basic features and types of computations in a simple
way, if not finer aspects which can also be considered.
\end{abstract}

	Let us consider first curves in the plane ${\bf R}^2$.  We
shall identify ${\bf R}^2$ with the set ${\bf C}$ of complex numbers.

	Let $\Gamma$ be some kind of curve in ${\bf C}$, or perhaps
union of pieces of curves.  For each $z \in {\bf C} \backslash \Gamma$,
we have the contour integral
\begin{equation}
\label{int_{Gamma} frac{1}{(z-zeta)^2} d zeta}
	\int_{\Gamma} \frac{1}{(z-\zeta)^2} \, d\zeta
\end{equation}
as from complex analysis.  More precisely, ``$d\zeta$'' is the element
of integration such that if $\gamma$ is an arc in ${\bf C}$ from a
point $a$ to another point $b$, then
\begin{equation}
	\int_{\gamma} d\zeta = b - a.
\end{equation}
This works no matter how $\gamma$ goes from $a$ to $b$.  This is different
from integrating arclength, for which the element of integration is often
written $|d\zeta|$.  For this we have that
\begin{equation}
	\int_{\gamma} |d\zeta| = {\rm length}(\gamma),
\end{equation}
and this very much depends on the way that $\gamma$ goes from $a$ to
$b$.

	If $\Gamma$ a union of closed curves, then
\begin{equation}
\label{int_{Gamma} frac{1}{(z-zeta)^2} d zeta = 0}
	\int_{\Gamma} \frac{1}{(z-\zeta)^2} \, d\zeta = 0.
\end{equation}
This is a standard formula from complex analysis (an instance of
``Cauchy formulae''), and one can look at it in the following manner.
As a function of $\zeta$, $1/(z-\zeta)^2$ is the complex derivative in
$\zeta$ of $1/(z-\zeta)$,
\begin{equation}
	\frac{d}{d\zeta} \biggl(\frac{1}{z-\zeta}\biggr)
				= \frac{1}{(z-\zeta)^2}.
\end{equation}
If $\gamma$ is a curve from points $a$ to $b$ again, which does not pass
through $z$, then
\begin{equation}
\label{int_{gamma} frac{1}{(z-zeta)^2} d zeta = frac{1}{z-b} - frac{1}{z-a}}
	\int_{\gamma} \frac{1}{(z-\zeta)^2} \, d\zeta
		= \frac{1}{z-b} - \frac{1}{z-a}.
\end{equation}
In particular, one gets $0$ for closed curves (since that corresponds to
having $a = b$).

	As a variation of these matters, if $\Gamma$ is a line, then
\begin{equation}
\label{int_{Gamma} frac{1}{(z-zeta)^2} d zeta = 0, 2}
	\int_{\Gamma} \frac{1}{(z-\zeta)^2} \, d\zeta = 0
\end{equation}
again.  This can be derived from (\ref{int_{gamma} frac{1}{(z-zeta)^2}
d zeta = frac{1}{z-b} - frac{1}{z-a}}) (and can be looked at in terms
of ordinary calculus, without complex analysis).  There is enough
decay in the integral so that there is no problem with using the whole
line.

	What would happen with these formulae if we replaced the
complex element of integration $d\zeta$ with the arclength element of
integration $|d\zeta|$?  In general we would not have
(\ref{int_{Gamma} frac{1}{(z-zeta)^2} d zeta = 0}) for unions of
closed curves, or (\ref{int_{gamma} frac{1}{(z-zeta)^2} d zeta =
frac{1}{z-b} - frac{1}{z-a}}) for a curve $\gamma$ from $a$ to $b$.
However, we would still have (\ref{int_{Gamma} frac{1}{(z-zeta)^2} d
zeta = 0, 2}) for a line, because in this case $d\zeta$ would be a
constant times $|d\zeta|$.

	Let us be a bit more general and consider an element of
integration $d\alpha(\zeta)$ which is positive, like the arclength
element $|d\zeta|$, but which is allowed to have variable density.
Let us look at an integral of the form
\begin{equation}
\label{int_{Gamma} frac{1}{(z-zeta)^2} dalpha(zeta)}
	\int_{\Gamma} \frac{1}{(z-\zeta)^2} \, d\alpha(\zeta).
\end{equation}
This integral can be viewed as a kind of measurement of \emph{curvature}
of $\Gamma$ (which also takes into account the variability of the density
in $d\alpha(\zeta)$).

	If we put absolute values inside the integral, then the result
would be roughly $\dist(z,\Gamma)^{-1}$,
\begin{equation}
\label{int_{Gamma} frac{1}{|z-zeta|^2} dalpha(zeta) approx dist(z,Gamma)^{-1}}
	\int_{\Gamma} \frac{1}{|z-\zeta|^2} \, d\alpha(\zeta)
			\approx \dist(z,\Gamma)^{-1}
\end{equation}
under suitable conditions on $\Gamma$.  For instance, if $\Gamma$ is a
line, then the left side of (\ref{int_{Gamma} frac{1}{|z-zeta|^2}
dalpha(zeta) approx dist(z,Gamma)^{-1}}) is equal to a positive
constant times the right side of (\ref{int_{Gamma} frac{1}{|z-zeta|^2}
dalpha(zeta) approx dist(z,Gamma)^{-1}}).  

	The curvature of a curve is defined in terms of the derivative
of the unit normal vector along the curve, or, what is essentially the
same here, the derivative of the unit tangent vector.  The unit tangent
vector gives exactly what is missing from $|d\zeta|$ to get $d\zeta$,
if we write the unit tangent vector as a complex number.  (One should
also follow the tangent in the orientation of the curve.)

	If the curve is a line or a line segment, then the tangent
is constant, which one can pull in and out of the integral.  In general
one can view (\ref{int_{Gamma} frac{1}{(z-zeta)^2} dalpha(zeta)}) as
a measurement of the variability of the unit tangent vectors, and of
the variability of the positive density involved in $d\alpha(\zeta)$.

	Let us look at some simple examples.  Suppose first that $\Gamma$
is a straight line segment from a point $a \in {\bf C}$ to another point
$b$, $a \ne b$.  Then $|d\zeta|$ is a constant multiple of $d\zeta$, and
\begin{equation}
	\int_{\Gamma} \frac{1}{(z-\zeta)^2} \, |d\zeta|
	   = ({\rm Constant}) \cdot \Bigl(\frac{1}{z-b} - \frac{1}{z-a}\Bigr).
\end{equation}
In this case the ordinary curvature is $0$, except that one can say
that there are contributions at the endpoints, like Direc delta
functions, which are reflected in right side.  If $z$ gets close to
$\Gamma$, but does not get close to the endpoints $a$, $b$ of
$\Gamma$, then the right side stays bounded and behaves nicely.  This
is ``small'' in comparison with $\dist(z,\Gamma)^{-1}$.  Near $a$ or
$b$, we get something which is indeed like $|z-a|^{-1}$ or
$|z-b|^{-1}$.

	As another example, suppose that we have a third point $p \in
{\bf C}$, where $p$ does not lie in the line segment between $a$ and
$b$ (and is not equal to $a$ or $b$).  Consider the curve $\Gamma$
which goes from $a$ to $p$ along the line segment between them, and
then goes from $p$ to $b$ along the line segment between them.  Again
$|d\zeta|$ is a constant multiple of $d\zeta$.  Now we have
\begin{equation}
	\int_{\Gamma} \frac{1}{(z-\zeta)^2} \, |d\zeta|
	   = c_1 \Bigl(\frac{1}{z-p} - \frac{1}{z-a}\Bigr)
		+ c_2 \Bigl(\frac{1}{z-b} - \frac{1}{z-p}\Bigr),
\end{equation}
where $c_1$ and $c_2$ are constants which are not equal to each other.
This is like the previous case, except that the right side behaves like
a constant times $|z-p|^{-1}$ near $p$ (and remains bounded away from
$a$, $b$, $p$).  This reflects the presence of another Dirac delta
function for the curvature, at $p$.  If the curve flattens out, so that
the angle between the two segments is close to $\pi$, then the coefficient
$c_1 - c_2$ of the $(z-p)^{-1}$ term becomes small.

	Now suppose that $\Gamma$ is the unit circle in ${\bf C}$, centered
around the origin.  In this case $|d\zeta|$ is the same as $d\zeta/\zeta$,
except for a constant factor, and we consider
\begin{equation}
	\int_{\Gamma} \frac{1}{(z-\zeta)^2} \, \frac{d\zeta}{\zeta}.
\end{equation}
If $z = 0$, then one can check that this integral is $0$.  For $z \ne 0$,
let us rewrite the integral as
\begin{equation}
	\int_{\Gamma} \frac{1}{(z-\zeta)^2} \, 
			\Bigl(\frac{1}{\zeta} - \frac{1}{z}\Bigr) \, d\zeta
	+ \frac{1}{z} \int_{\Gamma} \frac{1}{(z-\zeta)^2} \, d\zeta.
\end{equation}
The second integral is $0$ for all $z \in {\bf C} \backslash \Gamma$,
as in the earlier discussion.  The first integral is equal to
\begin{equation}
   \int_{\Gamma} \frac{1}{(z-\zeta)^2} \, \frac{(z-\zeta)}{z \zeta} \, d\zeta
= \frac{1}{z} \int_{\Gamma} \frac{1}{(z-\zeta)} \, \frac{1}{\zeta} \, d\zeta.
\end{equation}
On the other hand,
\begin{equation}
	\frac{1}{(z-\zeta)} \, \frac{1}{\zeta} 
	    = \frac{1}{z} \Bigl(\frac{1}{z-\zeta} + \frac{1}{\zeta}\Bigr),
\end{equation}
and so we obtain 
\begin{equation}
   \frac{1}{z^2} \int_{\Gamma} \Bigl(\frac{1}{z-\zeta} + \frac{1}{\zeta}\Bigr)
								 \, d\zeta.
\end{equation}
For $|z| > 1$ we have that
\begin{equation}
	\int_{\Gamma} \frac{1}{z-\zeta} \, d\zeta = 0,
\end{equation}
and thus we get a constant times $1/z^2$ above.  If $|z| < 1$, then
\begin{equation}
	\int_{\Gamma} \frac{1}{z-\zeta} \, d\zeta 
		= - \int_{\Gamma} \frac{1}{\zeta} \, d\zeta,
\end{equation}
and the earlier expression is equal to $0$.

	For another example, fix a point $q \in {\bf C}$, and suppose that
$\Gamma$ consists of a finite number of rays emanating from $q$.  On
each ray, we assume that we have an element of integration $d\alpha(\zeta)$
which is a positive constant times the arclength element.  

	If $R$ is one of these rays, then 
\begin{equation}
	\int_R \frac{1}{(z-\zeta)^2} \, d\alpha(\zeta)
		= ({\rm constant}) \, \frac{1}{z-q}.
\end{equation}
This constant takes into account both the direction of the ray and the
density factor in $d\alpha(\zeta)$ on $R$.

	If we now sum over the rays, we still get 
\begin{equation}
	\int_{\Gamma} \frac{1}{(z-\zeta)^2} \, d\alpha(\zeta)
		= ({\rm constant}) \, \frac{1}{z-q};
\end{equation}
however, this constant can be $0$.  This happens if $\Gamma$ is a
union of lines through $q$, with constant density on each line, and it
also happens more generally, when the directions of the rays satisfy a
suitable balancing condition, depending also on the density factors
for the individual rays.  This can happen with $3$ rays, for instance.

	When the constant is $0$, $\Gamma$ (with these choices of
density factors) has ``curvature $0$'', even if this is somewhat
complicated, because of the singularity at $q$.  This is a special
case of the situation treated in \cite{AA}.

	In general, ``weak'' or integrated curvature is defined using
suitable test functions on ${\bf R}^2$ with values in ${\bf R}^2$ (or
on ${\bf R}^n$ with values in ${\bf R}^n$), as in \cite{AA}.  For $n =
2$ one can reformulate this (trivially) in terms of complex-valued
functions on ${\bf C}$, and complex-analyticity gives rise to simpler
formulas.  This is what happens here.

	For more information on these topics, see Chapter 2 of Part
III of \cite{DS2}.  In \cite{DS2} there are further issues which are not
needed in various settings.
	
	Now let us look at similar matters in ${\bf R}^n$, $n > 2$,
and $(n-1)$-dimensional surfaces there.  Ordinary complex analysis is
no longer available, but there are substitutes, in terms of
quarternions (in low dimensions) and Clifford algebras.  For the sake
of definiteness let us focus on the latter.

	Let $n$ be a positive integer.  The \emph{Clifford algebra}
$\mathcal{C}(n)$ has $n$ generators $e_1, e_2, \ldots, e_n$ which
satisfy the following relations:
\begin{eqnarray}
	e_j \, e_k & = & - e_k \, e_j \quad\hbox{when } j \ne k;	\\
	e_j^2 & = & -1	       \qquad\hbox{ for all } j.	\nonumber
\end{eqnarray}
Here $1$ denotes the identity element in the algebra.  These are the
only relations.  More precisely, one can think of $\mathcal{C}(n)$ first
as a real vector space of dimension $2^n$, in which one has a basis
consisting of all products of $e_j$'s of the form
\begin{equation}
	e_{j_1} \, e_{j_2} \cdots e_{j_\ell},
\end{equation}
where $j_1 < j_2 < \cdots < j_\ell$, and $\ell$ is allowed to range from
$0$ to $n$, inclusively.  When $\ell = 0$ this is interpreted as
giving the identity element $1$.  If $\beta, \gamma \in
\mathcal{C}(n)$, then $\beta$ and $\gamma$ are given by linear
combinations of these basis elements, and it is easy to define the
product $\beta \, \gamma$ using the relations above and standard rules
(associativity and distributivity).

	If $n = 1$, then the result is isomorphic to the complex
numbers in a natural way, and if $n = 2$, the result is isomorphic to
the quarternions.  Note that $\mathcal{C}(n)$ contains ${\bf R}$
in a natural way, as multiples of the identity element.

	A basic feature of the Clifford algebra $\mathcal{C}(n)$ is
that if $\beta \in \mathcal{C}(n)$ is in the linear span of $e_1, e_2,
\ldots, e_n$ (without taking products of the $e_j$'s), then $\beta$
can be inverted in the algebra if and only if $\beta \ne 0$.  More
precisely, if
\begin{equation}
	\beta = \sum_{j=1}^n \beta_j \, e_j,
\end{equation}
where each $\beta_j$ is a real number, then
\begin{equation}
	\beta^2 = - \sum_{j=1}^n |\beta_j|^2.
\end{equation}
If $\beta \ne 0$, then the right side is a nonzero real number,
and $-(\sum_{j=1}^n |\beta_j|^2)^{-1} \beta$ is the multiplicative
inverse of $\beta$.

	More generally, if $\beta$ is in the linear span of
$1$ and $e_1, e_2, \ldots, e_n$, so that
\begin{equation}
	\beta = \beta_0 + \sum_{j=1}^n \beta_j \, e_j,
\end{equation}
where $\beta_0, \beta_1, \ldots, \beta_n$, then we set
\begin{equation}
	\beta^* = \beta_0 - \sum_{j=1}^n \beta_j \, e_j.
\end{equation}
This is analogous to complex conjugation of complex numbers, and
we have that
\begin{equation}
	\beta \, \beta^* = \beta^* \, \beta = \sum_{j=0}^n |\beta_j|^2.
\end{equation}
If $\beta \ne 0$, then $(\sum_{j=0}^n |\beta_j|^2)^{-1} \beta^*$ is
the multiplicative inverse of $\beta$, just as in the case of complex
numbers.

	When $n > 2$, nonzero elements of $\mathcal{C}(n)$ may not be
invertible.  For real and complex numbers and quarternions it is true
that nonzero elements are invertible.  The preceding observations are
substitutes for this which are often sufficient.

	Now let us turn to \emph{Clifford analysis}, which is an
analogue of complex numbers in higher dimensions using Clifford
algebras.  (See \cite{BDS} for more information.)

	Suppose that $f$ is a function on ${\bf R}^n$, or some
subdomain of ${\bf R}^n$, which takes values in $\mathcal{C}(n)$.  We
assume that $f$ is smooth enough for the formulas that follow (with
the amount of smoothness perhaps depending on the circumstances).
Define a differential operator $\mathcal{D}$ by
\begin{equation}
	\mathcal{D} f = \sum_{j=1}^n e_j \, \frac{\partial}{\partial x_j} f.
\end{equation}

	Actually, there are some natural variants of this to also
consider.  This is the ``left'' version of the operator; there is also
a ``right'' version, in which the $e_j$'s are moved to the right side of
the derivatives of $f$.  This makes a difference, because the Clifford
algebra is not commutative, but the ``right'' version enjoys the same
kind of properties as the ``left'' version.  (Sometimes one uses the
two at the same time, as in certain integral formulas, in which the
two operators are acting on separate functions which are then part of
the same expression.)  

	As another alternative, one can use the Clifford algebra
$\mathcal{C}(n-1)$ for Clifford analysis on ${\bf R}^n$, with one
direction in ${\bf R}^n$ associated to the multiplicative identity
element $1$, and the remaining $n-1$ directions associated to $e_1,
e_2, \ldots, e_{n-1}$.  There is an operator analogous to
$\mathcal{D}$, and properties similar to the ones that we are about to
describe (with adjustments analogous to the conjugation operation
$\beta \mapsto \beta^*$).

	For the sake of definiteness, let us stick to the version that
we have.  A function $f$ as above is said to be \emph{Clifford analytic}
if
\begin{equation}
	\mathcal{D} f = 0
\end{equation}
(on the domain of $f$).

	Clifford analytic functions have a lot of features analogous
to those of complex analytic functions, including integral formulas.
There is a natural version of a \emph{Cauchy kernel}, which is given
by
\begin{equation}
	\mathcal{E}(x-y) = \frac{\sum_{j=1}^n (x_j - y_j) \, e_j}{|x-y|^n}.
\end{equation}
This function is Clifford analytic in $x$ and $y$ away from $x = y$,
and it has a ``fundamental singularity'' at $x = y$, just as $1/(z-w)$
has in the complex case.

	One can calculate these properties directly, and one can also
look at them in the following way.  A basic indentity involving
$\mathcal{D}$ is 
\begin{equation}
	\mathcal{D}^2 = - \Delta,
\end{equation}
where $\Delta$ denotes the Laplacian, $\Delta = \sum_{j=1}^n \partial^2 /
\partial x_j^2$.  The kernel $\mathcal{E}(x)$ is a constant multiple of
\begin{eqnarray}
	&& \mathcal{D}(|x|^{n-2}) \quad\hbox{when } n > 2,		\\
	&& \mathcal{D}(\log |x|)  \quad\hbox{when } n = 2.	\nonumber
\end{eqnarray}
For instance, the Clifford analyticity of $\mathcal{E}(x)$ for $x \ne 0$
follows from the harmonicity of $|x|^{n-2}$, $\log |x|$ for $x \ne 0$
(when $n > 2$, $n = 2$, respectively).

	Analogous to (\ref{int_{Gamma} frac{1}{(z-zeta)^2} d zeta}), let
us consider integrals of the form
\begin{equation}
\label{int_{Gamma} frac{partial}{partial x_m} mathcal{E}(x-y) N(y) dy}
	\int_{\Gamma} \frac{\partial}{\partial x_m} \mathcal{E}(x-y) \,
 		     N(y) \, dy,   \quad x \in {\bf R}^n \backslash \Gamma,
\end{equation}
where $\Gamma$ is some kind of $(n-1)$-dimensional surface in ${\bf R}^n$,
or union of pieces of surfaces, 
\begin{equation}
	N(y) = \sum_{j=1}^n N_j(y) \, e_j
\end{equation}
is the unit normal to $\Gamma$ (using some choice of orientation for
$\Gamma$), turned into an element of $\mathcal{C}(n)$ using the
$e_j$'s in this way, and $dy$ denotes the usual element of surface
integration on $\Gamma$.  Thus $N(y) \, dy$ is a
Clifford-algebra-valued element of integration on $\Gamma$ which is
analogous to $d\zeta$ for complex contour integrals, as in
(\ref{int_{Gamma} frac{1}{(z-zeta)^2} d zeta}).  A version of the
Cauchy integral formula implies that
\begin{equation}
	\int_{\Gamma} \mathcal{E}(x-y) \, N(y) \, dy
\end{equation}
is locally constant on ${\bf R}^n \backslash \Gamma$ when $\Gamma$ is
a ``closed surface'' in ${\bf R}^n$, i.e., the boundary of some
bounded domain (which is reasonably nice).  In fact, this integral is
a nonzero constant inside the domain, and it is zero outside the
domain.  At any rate, the differentiated integral (\ref{int_{Gamma}
frac{partial}{partial x_m} mathcal{E}(x-y) N(y) dy}) is then $0$ for
all $x \in {\bf R}^n \backslash \Gamma$, in analogy with
(\ref{int_{Gamma} frac{1}{(z-zeta)^2} d zeta = 0}).

	Now suppose that we have a positive element of integration
$d\alpha(y)$ on $\Gamma$, which is the usual element of surface
integration $dy$ together with a positive density which is allowed
to be variable.  Consider integrals of the form
\begin{equation}
\label{int_{Gamma} frac{partial}{partial x_m} mathcal{E}(x-y) d alpha(y)}
	\int_{\Gamma} \frac{\partial}{\partial x_m} \mathcal{E}(x-y) 
 		     \, d\alpha(y),   \quad x \in {\bf R}^n \backslash \Gamma.
\end{equation}
This again can be viewed in terms of integrations of curvatures of
$\Gamma$ (also incorporating the variability of the density in
$d\alpha(y)$).  In a ``flat'' situation, as when $\Gamma$ is an
$(n-1)$-dimensional plane, or a piece of one, $N(y)$ is constant, and
if $d\alpha(y)$ is replaced with a constant times $dy$, then we can
reduce to (\ref{int_{Gamma} frac{partial}{partial x_m} mathcal{E}(x-y)
N(y) dy}), where special integral formulas such as Cauchy formulas
can be used.

	Topics related to this are discussed in Chapter 3 of Part III
of \cite{DS2}, although, as before, further issues are involved there
which are not needed in various settings.  See \cite{BDS} for more
on Clifford analysis, including integral formulas.  Related matters
of curvature are investigated in \cite{H}.

	Information about some other connections between Cauchy
kernels and geometry can be found in the references.

\end{document}